\documentclass[11pt]{article}
\usepackage{latexsym}
\usepackage{amsmath}
\usepackage{amssymb}
\usepackage{amsthm}
\usepackage{epsfig}

\newtheorem{theorem}{Theorem}
\newtheorem{corollary}{Corollary}
\newtheorem{lemma}{Lemma}
\newtheorem{conjecture}{Conjecture}

\newtheorem{proposition}{Proposition}
\newtheorem{definition}{Definition}

\allowdisplaybreaks

\def\itc#1{{\em #1\/}}

\begin{document}

\title{Wide Partitions, Latin Tableaux, and Rota's Basis Conjecture}
\author{Timothy Y. Chow\protect\footnote{250 Whitwell Street \#2,
Quincy, MA 02169, \texttt{tchow@alum.mit.edu}}
\and C. Kenneth Fan\protect\footnote{\texttt{ckfan@alum.mit.edu}}
\and Michel X. Goemans\protect\footnote{Dept.\ of Mathematics,
Mass.\ Inst.\ of Technology, Cambridge, MA 02139,
\texttt{goemans@math.mit.edu}}
\and Jan Vondrak\protect\footnote{Dept.\ of Mathematics,
Mass.\ Inst.\ of Technology, Cambridge, MA 02139,
\texttt{vondrak@math.mit.edu}}
}
\date{September 9, 2002}

\maketitle

\section{Introduction}

The main purpose of this paper is to publicize,
and to present partial results on,
a new combinatorial conjecture 
of Brian Taylor and the first author.
We begin by stating the conjecture.
(We assume some knowledge of the terminology of integer partitions;
readers lacking this background should consult~\cite{Stanley_EC2}.)

\begin{definition}
An integer partition~$\mu$ is a \itc{subpartition} of
an integer partition~$\lambda$ (written \hbox{$\mu\subseteq\lambda$})
if the multiset of parts of~$\mu$ is a submultiset of
the multiset of parts of~$\lambda$.
Equivalently, the Young diagram of~$\mu$ is
obtained by deleting some rows from the Young diagram of~$\lambda$.
\end{definition}

\begin{definition}
An integer partition~$\lambda$ is \itc{wide} if $\mu\ge\mu'$
in dominance order for all $\mu\subseteq\lambda$.
Here $\mu'$ denotes the conjugate of~$\mu$.
\end{definition}

\begin{conjecture}[The Wide Partition Conjecture for Free Matroids]
An integer partition~$\lambda$ is wide if and only if
there exists a tableau of shape~$\lambda$ such that
(1)~for all~$i$, the entries in the $i$th row of the tableau
are precisely the integers from $1$ to~$\lambda_i$ inclusive, and
(2)~for all~$j$, the entries in the $j$th column of the tableau
are pairwise distinct.
\end{conjecture}

We believe that
the wide partition conjecture (or WPC for short) for free matroids
has intuitive appeal as stated.
However, the reader might prefer one of the following
equivalent formulations.

\begin{itemize}
\item
In the language of edge colorings,
it states that for bipartite graphs arising from wide partitions,
the set of all color-feasible sequences has a unique maximal element.
\item
In the language of network flows,
it states that certain integer multiflow problems
that are associated with wide partitions satisfy
a max-flow min-cut theorem
and have integral optimal solutions.
\item
In the language of the Greene-Kleitman theorem,
it states that the line graph of a bipartite graph
arising from a wide partition
has a stable set cover
that is simultaneously $k$-saturated for all~$k$.
\end{itemize}

More precise statements of these reformulations
will be given later.

As we explain later, the motivation for the WPC for free matroids
comes from Rota's basis conjecture, which in turn is motivated by
certain questions in classical invariant theory.
A curious consequence of this connection to invariant theory is that
the WPC for free matroids might actually be more interesting
if it is \itc{false} rather than true,
because then it would probably lead to new and unsuspected
identities in invariant theory.
We do not describe the invariant-theoretic connection
in detail in this paper, but hope to do so elsewhere.

Our main partial result is that the WPC for free matroids is true
for certain wide partitions
with only a small number of distinct part sizes.
We also show that certain graphs arising from wide partitions
satisfy a property called ``$\Delta$-conjugacy,''
which Greene and Kleitman famously showed was true
of comparability graphs.
This result seems interesting in its own right,
because graphs satisfying $\Delta$-conjugacy
are rather hard to come by~\cite{Fomin},
and our examples seem to be new.
Finally, we show that to prove the WPC for free matroids,
it suffices to consider self-conjugate shapes.

\section{Basic definitions}

We follow \cite{Stanley_EC2} for most of our notation
and terminology for (integer) partitions,
but the reader should note two important exceptions.
Firstly, the subpartition relation $\mu\subseteq\lambda$
defined above is different from the usual one.
Secondly, for us a \itc{tableau}
is simply a Young diagram with a positive integer entry in each cell;
there is no implicit condition of semi\-standardness.

Young diagrams may be identified with bipartite graphs in a natural way.
If $\lambda$ is a partition,
we define $G_\lambda$ to be the bipartite graph whose vertices
are the rows and columns of~$\lambda$ and that has an edge
between row~$i$ and column~$j$ if and only if $(i,j)$
is a cell of the Young diagram of~$\lambda$
(i.e., if and only if $j\le\lambda_i$).

Sometimes it is more convenient to consider
$L(G_\lambda)$, the \itc{line graph} of~$G_\lambda$,
than to consider $G_\lambda$ itself.
The vertices of~$L(G_\lambda)$ are
the cells of the Young diagram of~$\lambda$,
and two vertices are adjacent if the cells lie in the same row or column.

The Young diagram of~$\lambda$ may also be identified with
a 0-1 matrix with $\ell(\lambda)$ rows and $\lambda_1$ columns;
the $(i,j)$ entry is one if and only if $(i,j)$ is 
a cell of the Young diagram. 

We will employ all the above ways of looking at Young diagrams,
so the reader should get used to switching freely
between the different viewpoints.

There are two well-known theorems that we need later.
See~\cite{BvN, Gale, Ryser} for proofs.

\begin{proposition}[Gale-Ryser Theorem]
\label{pr:GR}
Let $\lambda$ be a partition of~$n$ with $\ell$ parts and
let $\mu$ be a partition of~$n$ with $m$ parts.
Then there exists an $\ell\times m$ 0-1 matrix~$A$
whose $i$th row sums to~$\lambda_i$ (for all~$i$) and
whose $j$th column sums to~$\mu_j$ (for all~$j$)
if and only if $\lambda' \ge \mu$.
\end{proposition}

\begin{proposition}[Birkhoff-von Neumann Theorem]
\label{pr:BvN}
A nonnegative integer square matrix whose rows and columns
all sum to~$n$ may be written as the sum of $n$ permutation matrices.
\end{proposition}

\section{Wide partitions}

As we said in the introduction,
a partition~$\lambda$ is \itc{wide} if $\mu\ge\mu'$ for all
$\mu\subseteq\lambda$.
In this section we prove some fundamental facts about wide partitions.

The number of wide partitions of~$n$ is an integer sequence that begins

1, 1, 2, 3, 3, 5, 6, 9, 11, 14, 18, 23, 29, 35, 45, 56, 68, 85, 103, 125,
150, 183, 217, 266, 315, 380, 449, 534, 628, 745, 874, 1034, 1212, 1423,
1665, 1944, 2265, 2627, 3055, 3536, 4099, 4735, 5479, 6309, 7273, 8358,
9599, 11012, 12605, 14421, 16480, 18825, 21456, 24474, 27822, 31677,
35934, 40825, 46217, 52420, 59253, 67056, 75699, 85532, 96407.

Superseeker does not recognize this sequence.

Ostensibly, checking wideness requires checking
all subpartitions, a potentially exponential-time computation.
We show next that checking wideness takes only polynomial time.

\begin{definition}
A subpartition $\mu\subseteq\lambda$ is a \itc{lower subpartition}
if $\mu$ is obtained from~$\lambda$ by deleting the largest~$i$
parts of~$\lambda$ for some~$i\ge 0$.
\end{definition}

The following fact was first conjectured by Xun Dong (personal
communication).

\begin{proposition}
\label{pr:xundong}
If $\lambda$ is a partition such that $\mu\ge\mu'$ for
all lower subpartitions $\mu$ of~$\lambda$, then $\lambda$ is wide.
\end{proposition}

\begin{proof}
If $\lambda$ is a partition, let $\lambda^i$ denote the
subpartition of~$\lambda$ obtained by deleting the $i$th part of~$\lambda$.
Thus $\lambda^i_j = \lambda_j$ if $j<i$ and $\lambda^i_j = \lambda_{j+1}$
if $j\ge i$.

The proof is by induction on the number of parts of~$\lambda$.
Let $\lambda$ be a partition such that $\mu\ge\mu'$ for
all lower subpartitions $\mu$ of~$\lambda$.
Then in particular, $\lambda\ge\lambda'$ and $\lambda^1\ge(\lambda^1)'$.
We claim that $\lambda^i \ge (\lambda^i)'$ for all $i$.
To see this, fix any~$i$.
We need to show that for all~$j$,
$$\sum_{k=1}^j \lambda^i_k \ge \sum_{k=1}^j (\lambda^i)'_k.$$
Note that it suffices to consider only those $j\le \lambda_1$,
so we henceforth assume that $j\le \lambda_1$.

If $j<i$ then because $\lambda\ge\lambda'$, we have
$$\sum_{k=1}^j \lambda^i_k
 = \sum_{k=1}^j \lambda_k
 \ge \sum_{k=1}^j \lambda'_k
\ge \sum_{k=1}^j (\lambda^i)'_k,$$
so let us suppose that $j\ge i$.  We split into two cases,
the first case being the case in which $j \le \lambda_i$.
Then
\begin{eqnarray*}
\sum_{k=1}^j \lambda^i_k
 &=&  \sum_{k=2}^{j+1} \lambda_k + (\lambda_1 - \lambda_i)\\
 &\ge& \sum_{k=2}^{j+1} \lambda_k\\
 &\ge& \sum_{k=1}^{j} (\lambda'_k-1)
   \qquad \mbox{(because $\lambda^1\ge (\lambda^1)'$ and
$j\le\lambda_1$)}\\
 &=& \sum_{k=1}^j (\lambda^i)'_k
   \qquad \mbox{(because $j\le\lambda_i$)}.
\end{eqnarray*}
In the second case, $j > \lambda_i$, so
\begin{eqnarray*}
\sum_{k=1}^j \lambda^i_k
 &=& \sum_{k=2}^{j+1} \lambda_k + (\lambda_1 - \lambda_i)\\
 &\ge& \sum_{k=2}^{j+1} \lambda_k + (j - \lambda_i)
   \qquad \mbox{(because $j\le\lambda_1$)}\\
 &\ge& \sum_{k=1}^{j} (\lambda'_k-1) + (j - \lambda_i)\\
 &=& \sum_{k=1}^{\lambda_i} (\lambda^i)'_k +
    \sum_{k=\lambda_i + 1}^j \biggl((\lambda^i)'_k - 1\biggr) +
    (j - \lambda_i)\\
 &=& \sum_{k=1}^j (\lambda^i)'_k.
\end{eqnarray*}
This proves the claim.
Now note that by induction, $\lambda^1$ is wide.
It follows that $\lambda^i$ is wide for all~$i$,
because we have just shown that $\lambda^i \ge (\lambda^i)'$,
and every \itc{proper} lower subpartition~$\mu$ of~$\lambda^i$
is a subpartition of~$\lambda^1$ and therefore satisfies $\mu\ge\mu'$,
so we can again apply induction to conclude that $\lambda^i$ is wide.

Finally, suppose $\mu$ is a subpartition of~$\lambda$.
If $\mu=\lambda$ then $\mu\ge\mu'$ because $\lambda\ge\lambda'$.
Otherwise, $\mu\subseteq\lambda^i$ for some~$i$,
and therefore satisfies $\mu\ge\mu'$ because $\lambda^i$ is wide.
\end{proof}

The following easy but useful lemma has been independently observed
by several people, including D.~Waugh.

\begin{lemma}
If $\lambda$ is wide then
$\lambda_{\ell(\lambda)-i} > i$ for all $i\ge0$.
\end{lemma}

\begin{proof}
Since $\lambda$ is wide, so is the subpartition~$\mu$
consisting of the last $i+1$ rows of~$\lambda$.
The largest part of~$\mu$ is 
$\lambda_{\ell(\lambda)-i}$.
The first column of~$\mu$ is $i+1$.
Since $\mu\ge\mu'$, it follows that
$\lambda_{\ell(\lambda)-i} \ge i+1>i$.
\end{proof}

\begin{definition}
If $\lambda$ and~$\mu$ are partitions
then $\lambda+\mu$ denotes the partition whose $i$th part is
$\lambda_i+\mu_i$.
\end{definition}

\begin{proposition}
\label{pr:ckfan}
If $\lambda$ is wide and $\mu$ is a single column whose height is
at most $\lambda'_1+1$, then $\lambda+\mu$ is wide.
\end{proposition}

\begin{proof}
We claim that it suffices to show the following statement.

\begin{quote}
If $\lambda$ is wide and $\mu$ is a single column whose height is
at most $\lambda'_1+1$, then $\lambda+\mu \ge (\lambda+\mu)'$.
\end{quote}

For if we can prove this, then we can apply it to any subpartition
of our original partition~$\lambda$ to deduce the proposition.

Fix $i$.  We want to show that
the sum of the first $i$ rows of $\lambda+\mu$
is at least
the sum of the first $i$ columns of $\lambda+\mu$.
We split into two cases. 

{\em Case 1: $i\le\mu'_1$.}
In passing from $\lambda$ to $\lambda+\mu$,
the sum of the first $i$ rows increases by~$i$.
As for the columns, note that 
in passing from $\lambda$ to $\lambda+\mu$,
all we are doing is adding a column of height $\mu'_1$.
Therefore this causes the sum of the first $i$ columns
to increase by at most $\mu'_1-\lambda'_i$.
But by Lemma~1, $\lambda'_i\ge \lambda'_1-i+1$.
Therefore the increase in the sum of the first $i$ columns is at most
$$\mu'_1-\lambda'_i \le (\lambda'_1+1) - (\lambda'_1-i+1) = i,$$
which completes the proof of this case.

{\em Case 2: $i>\mu'_1$.}
In passing from $\lambda$ to $\lambda+\mu$,
the sum of the first $i$ rows increases by $\mu'_1$.
But the sum of the first $i$ columns cannot increase by more than~$\mu'_1$
either, so this case is also settled.
\end{proof}

\begin{corollary}
\label{co:sum}
If $\lambda$ and $\mu$ are wide then so is $\lambda+\mu$.
\end{corollary}

\begin{proof}
Since $\lambda+\mu=\mu+\lambda$ we may assume that $\lambda'_1\ge\mu'_1$.
Add the columns of~$\mu$ to $\lambda$ one by one, applying
Proposition~\ref{pr:ckfan} each time.
\end{proof}

\begin{definition}
A wide partition $\lambda$ is \itc{decomposable}
if there exist wide partitions $\mu$ and~$\nu$ such that
$\lambda=\mu+\nu$; it is \itc{indecomposable} otherwise.
\end{definition}

{\bf Caution.}
Although every wide partition is a sum of indecomposables,
the decomposition need not be unique.

\begin{proposition}
\label{pr:indec}
For any fixed $\ell$, the number of indecomposable wide partitions
with $\ell$ parts is finite.
\end{proposition}

Our proof of Proposition~\ref{pr:indec} uses the following lemma.

\begin{lemma}
\label{le:wide}
Let $\lambda$ be a wide partition of~$n$ and let $a$ be a positive integer.
Then for all sufficiently large~$b$, the partition
$$ \mu = (\overbrace{b, b, \ldots, b,}^{a\;\rm times} \lambda_1, \lambda_2,
  \ldots, \lambda_{\ell(\lambda)})$$
is wide.
\end{lemma}

\begin{proof}
We may obtain a weaker claim than Lemma~\ref{le:wide} by replacing
``$\mu$ is wide'' by the weaker conclusion ``$\mu\ge\mu'$.''
Proving this weaker claim suffices to prove the lemma,
because by Proposition~\ref{pr:xundong}
one need only check lower subpartitions of~$\mu$,
and all such lower subpartitions are either covered by the weaker claim
or are subpartitions of the wide partition~$\lambda$.

We now prove the weaker claim.
Write $\ell$ for $\ell(\mu)$.
Pick any $b \ge n/a + \ell$; we shall see that this is sufficiently large.
We want to show that for all~$i\le \ell$,
the sum of the first $i$ rows of~$\mu$
is at least
the sum of the first $i$ columns of~$\mu$.
We split into two cases.

{\em Case 1: $i\le a$.}
The sum of the first $i$ rows of~$\mu$ is $ib$.
The sum of the first $i$ columns of~$\mu$ is at most $i\ell$.
But $b \ge \ell$ by construction.

{\em Case 2: $a < i\le \ell$.}
The sum of the first $i$ rows of~$\mu$ is at least $ab$.
By choice of~$b$, $ab \ge n+a\ell$.
But $n+a\ell$ is at least
the sum of the first $\ell$ columns of~$\mu$
(since $n$ is large enough to encompass all of~$\lambda$,
and $a\ell$ is large enough to encompass the sum of the first
$\ell$ columns of the first $a$ rows of~$\mu$),
which in turn is at least
the sum of the first $i$ columns of~$\mu$,
since $\ell \ge i$.
\end{proof}

\begin{proof}[Proof of Proposition~\ref{pr:indec}]
Call a partition~$\mu$ \itc{squarish} if $\mu_{\ell(\mu)} \ge \ell(\mu)$.
Any squarish partition with $\ell$ parts 
may be obtained by starting with an $\ell \times \ell$ square shape
and adding columns of height at most~$\ell$ to it.
Therefore, by Proposition~\ref{pr:ckfan}, all squarish partitions are wide.

Let $\lambda$ be an indecomposable wide partition with $\ell$ parts.
We show by induction on~$i$ that
$\lambda_{\ell-i} - \lambda_{\ell-i+1}$ is bounded for all~$i\le\ell-1$.
This implies the proposition.

If $i=0$, then $\lambda_\ell \le 2\ell-1$; otherwise we would have
$\lambda = \mu+\nu$ with $\mu$ an $\ell\times\ell$ square and
$\nu$ a squarish partition.

For larger~$i$,
we know by induction that the lower subpartition~$\kappa$
consisting of the last $i$ parts of~$\lambda$
can only be one of a finite set of possible partitions.
For any fixed~$\kappa$,
observe that if $\lambda_{\ell-i} - \lambda_{\ell-i+1}$
is sufficiently large,
then we may write $\lambda=\mu+\nu$
where $\nu$ is a squarish partition with $\ell-i$ parts
and $\mu$ is of the form given in Lemma~\ref{le:wide}
(with the ``$\lambda$'' of Lemma~\ref{le:wide} being $\kappa$
and ``$a$'' being $\ell-i$).
So since $\lambda$ is an indecomposable wide partition,
$\lambda_{\ell-i} - \lambda_{\ell-i+1}$ is bounded.
There are only finitely many choices for~$\kappa$,
so the proof is complete.
\end{proof}

\section{Latin tableaux and the Wide Partition Conjecture}

\begin{definition}
If $M$ is a matroid, then an \itc{$M$-tableau} is a Young diagram with an
element of~$M$ in each cell of the diagram.
\end{definition}

\begin{definition}
Let $\lambda$ be a partition.
We say that $\lambda$ \itc{satisfies Rota's conjecture}
if, for any matroid~$M$ and any sequence $(I_i)$
of independent sets of~$M$ satisfying $|I_i| = \lambda_i$ for all~$i$,
there exists an $M$-tableau~$T$ of shape~$\lambda$ such that
(a)~for all~$i$,
the set of elements in the $i$th row of~$T$ is~$I_i$, and
(b)~for all~$j$,
the elements in the $j$th column of~$T$ comprise an independent set of~$M$.
(In particular, the elements in the $j$th column are pairwise distinct.)
\end{definition}

\begin{conjecture}[The Wide Partition Conjecture]
A partition~$\lambda$ satisfies Rota's conjecture if and only if it is wide.
\end{conjecture}

We shall see shortly that wideness is necessary; it is sufficiency that is
the real question.
The WPC contains Rota's basis conjecture~\cite{HR} as a special case.
It was formulated by Brian Taylor and the first author,
originally with the hope that it would allow Rota's basis conjecture
to be proved by induction on the number of cells in a wide partition.

Unfortunately, the WPC does not seem to be any easier than Rota's
basis conjecture.  Nevertheless, we believe that the WPC is interesting
in its own right, because
in the invariant-theoretic context that
originally motivated Rota's basis conjecture, there is nothing special
about square shapes.
If you believe Rota's basis conjecture,
then you should probably believe the WPC too.

Since the WPC seems hard, we have focused on
the special case of free matroids.

\begin{definition}
Let $\lambda$ and $\mu$ be partitions.
A \itc{Latin tableau~$T$ of shape~$\lambda$ and content~$\mu$}
is a Young diagram of shape~$\lambda$ with a single positive integer
in each cell such that (a)~no two cells in the same row or column
have the same entry and (b)~the total number of occurrences of
the integer~$i$ equals~$\mu_i$.
A partition~$\lambda$ is \itc{Latin} if there exists
a Latin tableau~$T$ of shape~$\lambda$ and content~$\lambda'$.
\end{definition}

It is not hard to see that in
a Latin tableau~$T$ of shape~$\lambda$ and content~$\lambda'$,
the entries in row~$i$ are precisely the integers from 
1 to~$\lambda_i$.  It follows that if $\lambda=\lambda'$, then in
a Latin tableau~$T$ of shape~$\lambda$ and content~$\lambda'=\lambda$,
the entries in \itc{column~$i$} are also precisely the integers from
1 to~$\lambda_i$.

{\bf The WPC for Free Matroids.} 
A partition~$\lambda$ is Latin if and only if it is wide.

We have verified the WPC for free matroids
by computer for all partitions whose Young diagram has at most 65 cells.
This set of partitions includes all indecomposable wide partitions
with at most five parts.
We have also verified the WPC for free matroids
for all partitions whose Young diagram
fits inside a $10\times10$ square.

Readers familiar with the Alon-Tarsi conjecture on Latin squares
may wonder if the WPC for matroids representable over a field
of characteristic zero follows from an Alon-Tarsi-like conjecture
that the number of ``even'' Latin tableaux is not equal to the number
of ``odd'' Latin tableaux of the same shape.
We expect this to be true and
provable in the same way that it is proved for square shapes,
but we have not verified the details.

As Victor Reiner was the first to observe, it is easy to see that
if $\lambda$ is Latin, then it is wide.
For let $T$ be a Latin tableau of shape~$\lambda$ and content~$\lambda'$.
If $\mu\subseteq\lambda$, then $T$ restricted to~$\mu$
is a Latin tableau---call it $S$---of shape~$\mu$ and content~$\mu'$.
We want to show that $\mu\ge\mu'$.
Fix $i$ and erase all the entries of~$S$ except those that are
less than or equal to~$i$.
There are at most $i$ entries remaining in each column of~$S$,
so if we ``push them up'' as far as possible,
we can fit them all into the first $i$ rows of~$S$.
Therefore the first $i$ rows of~$\mu$ contain at least
as many cells as the sum of the first $i$ parts of~$\mu'$.

As an aside, we remark that Latin tableaux,
especially of self-conjugate shapes,
seem to be quite pleasing structures.
Many concepts associated with Latin squares,
such as orthogonality,
can be generalized to Latin tableaux.
We speculate that Latin tableaux may have applications
to error-correcting codes and/or the design of experiments.

\section{Relationship with list coloring}

There is an alternative form of the WPC for free matroids,
which we now describe.

\begin{definition}
A partition $\lambda$ is \itc{strongly Latin}
if, for any sequence $(I_i)$ of sets of distinct integers
satisfying $|I_i| = \lambda_i$ for all~$i$,
there exists a tableau~$T$ of shape~$\lambda$ such that
(a)~for all~$i$,
the set of integers in the $i$th row of~$T$ is~$I_i$, and
(b)~for all~$j$,
the integers in the $j$th column of~$T$ are pairwise distinct.
\end{definition}

{\bf The WPC for Free Matroids, alternative form.}
A partition~$\lambda$ is strongly Latin if and only if it is wide.

If we recall the statement of the (full) WPC,
then this alternative form of the WPC for free matroids
might seem more natural than the form we stated in the previous section.
It matters little, however, since we shall see that the two forms are
equivalent.

It is clear that a strongly Latin partition is Latin.
One might think that the converse would be easy to prove
since intuitively the ``worst case'' is the one in which
the sets $I_i$ intersect as much as possible.
However, this is the same intuition that leads to
the false conclusion that the list chromatic number of a graph
must always equal its ordinary chromatic number.
Therefore we must tread carefully.

\begin{definition}
An \itc{orientation} of a graph~$G$ is an
assignment of a direction to each of the edges of~$G$.
\end{definition}

\begin{proposition}[Galvin]
\label{pr:galvin}
Let $G$ be the line graph of a bipartite graph,
and suppose that
each vertex of~$G$ is equipped with a list of available colors.
If there exists an orientation of~$G$
in which every complete subgraph of~$G$ is acyclic and
in which the outdegree of every vertex
is less than the number of (distinct) colors in its list,
then there is a list coloring of~$G$
(i.e., a choice, for each vertex, of a color from its list
in such a way that distinct colors are chosen for adjacent vertices).
\end{proposition}

\begin{proof}
See \cite{Galvin}.
\end{proof}

\begin{theorem}
\label{th:strong}
If $\lambda$ is Latin then it is strongly Latin.
\end{theorem}

\begin{proof}
Assume that $\lambda$ is Latin, so that there exists
a Latin tableau~$T$ of shape~$\lambda$ and content~$\lambda'$.
Use $T$ to define an orientation of~$L(G_\lambda)$, as follows:
Let an edge between two cells in the same \itc{row}
point to the cell whose entry in~$T$ is {\it larger,} and
let an edge between two cells in the same \itc{column}
point to the cell whose entry in~$T$ is {\it smaller.}
It is easily verified that in this orientation,
the outdegree of a vertex in the $i$th row is at most $\lambda_i-1$.

To see that $\lambda$ is strongly Latin,
suppose we are given a sequence~$(I_i)$ of sets of distinct integers
satisfying $|I_i| = \lambda_i$.
The existence of the tableau in the definition of ``strongly Latin''
is equivalent to the existence of a list coloring of~$L(G_\lambda)$
if each vertex in row~$i$ of~$L(G_\lambda)$ is equipped with the list~$I_i$.
So the orientation of~$L(G_\lambda)$ constructed above,
combined with Proposition~\ref{pr:galvin},
implies the theorem.
\end{proof}

Theorem~\ref{th:strong} becomes easier to prove
if we restrict ourselves to square shapes.
Two direct proofs of this special case were given in~\cite{Chow},
and it also follows immediately from
the Lebensold-Fulkerson theorem \cite{Fulk, Lebensold}
on disjoint matchings in bipartite graphs.

We remark that Galvin's theorem allows us to
prove something slightly stronger than Theorem~\ref{th:strong}.
Say that an orientation of~$L(G_\lambda)$ is \itc{colorable}
if every complete subgraph is acyclic and
the outdegree of a vertex in the $i$th row is at most $\lambda_i-1$.
Galvin tells us that to prove that $\lambda$ is strongly Latin,
we need only construct a colorable orientation.
This can be done using something slightly weaker than the Latin property.

\begin{definition}
A tableau~$T$ of shape~$\lambda$ is \itc{weakly Latin} if
(a)~for all~$i$,
the set of integers in the $i$th row of~$T$
is~$\{1, 2, \ldots, \lambda_i\}$, and
(b)~for all~$j$ and~$k$,
there are at most~$k$ entries in the $j$th column of~$T$
that are less than or equal to~$k$.
A partition $\lambda$ is \itc{weakly Latin}
if there exists a weakly Latin tableau of shape~$\lambda$.
\end{definition}

\begin{proposition}
\label{pr:weak}
A partition is weakly Latin if and only if it is Latin.
\end{proposition}

\begin{proof}
Essentially the same construction as above
shows that if $\lambda$ is weakly Latin then
there exists a colorable orientation of~$L(G_\lambda)$.
\end{proof}

We conclude this section with an application of the above results.

\begin{proposition}
\label{pr:sum}
If $\lambda$ and $\mu$ are Latin then so is $\lambda+\mu$.
\end{proposition}

\begin{proof}
Assume that $\lambda$ and $\mu$ are Latin.
Then by Theorem~\ref{th:strong}, $\mu$ is strongly Latin.
Let $T_\lambda$ be a Latin tableau of shape~$\lambda$ and
content~$\lambda'$.
Let $T_\mu$ be a Latin tableau whose $i$th row contains the integers
$\lambda_i+1, \lambda_i+2, \ldots, \lambda_i+\mu_i$ in some order.
Such a $T_\mu$ exists because $\mu$ is strongly Latin.
If we now take the union of
the set of columns of~$T_\lambda$ with
the set of columns of~$T_\mu$,
sort the columns according to height,
and combine them to form a tableau~$T$
of shape $\lambda+\mu$,
then we see that~$T$ is in fact a Latin tableau
of shape $\lambda+\mu$ and content $(\lambda+\mu)'$.
\end{proof}

\begin{corollary}
\label{co:indec}
If all indecomposable wide partitions with $\ell$ parts are Latin
then all wide partitions with $\ell$ parts are Latin.
\end{corollary}

\begin{proof}
This follows from Corollary~\ref{co:sum} and Proposition~\ref{pr:sum}.
\end{proof}

Our computer calculation therefore shows that all wide partitions
with at most five parts are Latin.
Unfortunately, the set of indecomposable wide partitions
does not seem to be any more tractable than
the set of all wide partitions,
so at this point it is not clear how useful Corollary~\ref{co:indec}~is.

\section{Relationship with the Greene-Kleitman theorem}

Much of what follows can be stated in the general framework of
antiblocking hypergraphs,
but for simplicity we restrict our attention to the case of
perfect graphs.
Readers unfamiliar with the terminology of perfect graphs
can find complete definitions in~\cite{Saks},
which we shall be citing several times.

Let $G$ be a perfect graph.
A \itc{$k$-clique} is a union of $k$ complete subgraphs of~$G$, and
a \itc{$k$-stable set} is a union of $k$ stable sets of~$G$.
We let $\omega_k(G)$ denote the cardinality (number of vertices)
of the largest $k$-clique of~$G$ and
we let $\alpha_k(G)$ denote the cardinality
of the largest $k$-stable set of~$G$.
We also define
$$\Delta\omega_k(G) = \omega_k - \omega_{k-1} \qquad\hbox{and}\qquad
  \Delta\alpha_k(G) = \alpha_k - \alpha_{k-1}.$$
If there is no confusion, then we may drop the~``$G$'' from the notation
for simplicity.

If $ \Delta\omega $ and $\Delta\alpha$ are partitions
(i.e., $\Delta\omega_1 \ge  \Delta\omega_2 \ge \Delta\omega_3 \ge \cdots$
and $\Delta\alpha_1 \ge  \Delta\alpha_2 \ge \Delta\alpha_3 \ge \cdots$)
and furthermore are conjugates of each other,
then we say that $G$ \itc{satisfies $\Delta$-conjugacy}.
It is a famous theorem, due to Greene and Kleitman \cite{Greene, GK},
that comparability graphs of finite posets satisfy $\Delta$-conjugacy.

A \itc{clique cover} of~$G$ is
a vertex-disjoint union of complete subgraphs
whose union covers all vertices of~$G$.
If $\lambda$ is a clique cover, then we abuse notation
and also let $\lambda$ denote the integer partition consisting of
the sizes of the cliques (arranged in nonincreasing order of course).
If $\lambda_k = \Delta\omega_k$ for all~$k$,
then we say that $\lambda$ is a \itc{uniform} clique cover.
(In general, uniform clique covers need not exist.)
We define \itc{stable set covers}
and \itc{uniform stable set covers}
in a completely analogous way.

Let $k$ be a positive integer.
A clique cover~$\lambda$ is \itc{$k$-saturated} if
$$\alpha_k = \sum_{i=1}^k \lambda'_i.$$
If $\lambda$ is simultaneously $k$-saturated for all~$k$,
then we say that $\lambda$ is \itc{completely saturated}.
Similarly a stable set cover~$\lambda$ is $k$-saturated if
$$\omega_k = \sum_{i=1}^k \lambda'_i,$$
and is completely saturated if it is $k$-saturated for all~$k$.
For arbitrary graphs, $k$-saturated clique/stable set covers
need not exist,
and even for comparability graphs,
completely saturated clique/stable set covers need not exist.

\begin{proposition}
\label{pr:Saks}
If $G$ is a perfect graph satisfying $\Delta$-conjugacy,
then for every positive integer $k$,
there exists a clique cover that is
simultaneously $k$-saturated  and $(k+1)$-saturated, and
there also exists a stable set cover that is
simultaneously $k$-saturated  and $(k+1)$-saturated.
\end{proposition}

\begin{proof}
Theorem 4.13 of \cite{Saks}.
\end{proof}

The conclusion of Proposition~\ref{pr:Saks}
is sometimes referred to as the \itc{t-phenomenon}.

The concept of uniform clique/stable set covers
does not seem to be as standard as the other concepts above.
We have not found a reference for the following simple proposition,
although it is unlikely to be new.

\begin{proposition}
\label{pr:uniform}
Let $G$ be a perfect graph.
Every completely saturated clique cover
is uniform.
If for all $k$ there exists
a $k$-saturated clique cover,
then every uniform clique cover
is completely saturated.
Both statements hold with ``stable set'' in place of ``clique.''
\end{proposition}

\begin{proof}
The complement of a perfect graph is perfect~\cite{Lovasz},
so it suffices to consider clique covers.

Let $\lambda$ be a completely saturated clique cover.
Fix~$k$.
There exists a $\lambda_k$-stable set~$S$
with cardinality $\sum_{i=1}^{\lambda_k} \lambda'_i$.
Now, $S$ contains at most $\min(\lambda_k, \lambda_i)$ vertices from
the $i$th clique of~$\lambda$.
But the cardinality of~$S$ forces $S$ to contain
\itc{exactly} $\min(\lambda_k, \lambda_i)$ vertices from
the $i$th clique of~$\lambda$.
Therefore, each of the $k$ largest cliques of~$\lambda$
(which all have cardinality at least~$\lambda_k$)
contains one element from each stable set of~$S$.
It follows that each stable set of~$S$ has at least $k$~vertices.

Now augment $S$ to a stable set \itc{cover~$S^+$} by adjoining
singleton sets.
These singletons are precisely the vertices in
the $k$ largest cliques of~$\lambda$ that are \itc{not} in~$S$.
Therefore, for any $k$-clique~$C$---in particular, one of maximum
cardinality---we have
$$|C| \le \sum_{s\in S^+} \min(k, |s|) 
  = \sum_{s\in S} \min(k,|s|) + \sum_{i=1}^k (\lambda_i - \lambda_k)
  = k\lambda_k  + \sum_{i=1}^k (\lambda_i - \lambda_k)
  = \sum_{i=1}^k \lambda_i.$$
Since $k$ was arbitrary, $\lambda$ is uniform.

Conversely, let $\lambda$ be a uniform clique cover.
Fix~$k$ and let $\mu$ be a $k$-saturated clique cover.
Because $\lambda$ is uniform, $\lambda\ge\mu$,
i.e., $\lambda'\le\mu'$, so in particular
$$ \sum_{i=1}^k \lambda'_i \le \sum_{i=1}^k \mu'_i.$$
Because $\mu$ is $k$-saturated,
there exists a $k$-stable set~$S$ such that
$$\sum_{i=1}^k \mu'_i = |S|.$$
Finally, because $\lambda$ is a clique cover,
$$ |S| \le \sum_{i=1}^k \lambda'_i.$$
Combining these facts forces the inequalities to be equalities,
and therefore $\lambda$ is $k$-saturated.
Since $k$ was arbitrary, $\lambda$ is completely saturated.
\end{proof}

Line graphs of bipartite graphs enjoy certain properties
that arbitrary perfect graphs do not,
as the following proposition illustrates.

\begin{proposition}
\label{pr:partition}
If $G$ is the line graph of a bipartite graph,
then $\Delta\alpha$ is a partition,
and for every positive integer~$k$,
there exists a $k$-saturated clique cover of~$G$.
Moreover, if $\Delta\omega$ is a partition,
then $G$ satisfies $\Delta$-conjugacy.
\end{proposition}

\begin{proof}
Theorems 4.18 and 4.23 of~\cite{Saks}.
(That $\Delta\alpha$ is a partition was already proved
in Lemma~2.1 of~\cite{deWerra}.)
\end{proof}

Not much beyond the conclusions of Proposition~\ref{pr:partition}
can be said, even if we require $G$ to equal $L(G_\lambda)$
for a (not necessarily wide) partition~$\lambda$.
For example,
if we take $\lambda = (7,7,6,6,3,3,3)$ and $G=L(G_\lambda)$,
then there is no uniform clique cover,
and in fact $\Delta\omega$ is not even a partition.
Moreover, there is no 5-saturated stable set cover.
However, one interesting question does remain open.

{\bf Latin Tableau Question.}
Let $G=L(G_\lambda)$ for an arbitrary partition~$\lambda$.
Does there necessarily exist a uniform stable set cover?

Note that line graphs of arbitrary bipartite graphs
need not have uniform stable set covers.
If the answer to the Latin Tableau Question is yes, then
this would not only verify the WPC for free matroids,
but would also give a necessary and sufficient condition for
the existence of a Latin tableau of shape~$\lambda$ and content~$\mu$,
for arbitrary $\lambda$ and~$\mu$.

If $\lambda$ is required to be wide,
then one easily deduces much stronger conclusions.

\begin{lemma}
\label{le:uniform}
If $\lambda$ is wide then the set of rows of
the Young diagram of~$\lambda$
is a uniform clique cover of~$L(G_\lambda)$.
\end{lemma}

\begin{proof}
It suffices to show that the maximum cardinality of any $k$-clique
is the sum of the first $k$ parts of~$\lambda$,
for all $k\le \ell(\lambda)$.
Let $C$ be a $k$-clique.
Since we are trying to maximize~$|C|$,
we may assume that the cliques of~$C$ are maximal.
Then $C$ is the union of $i$~rows and $j$~columns
for some nonnegative integers $i$ and~$j$ satisfying $i+j = k$.
Again, since we are trying to maximize~$|C|$,
we may assume that $C$ is the union of the \itc{first}
$i$~rows and the \itc{first} $j$~columns.
But because $\lambda$ is wide,
the lower subpartition~$\mu$ of~$\lambda$
comprising the last $\ell(\lambda)-i$ parts of~$\lambda$
satisfies $\mu \ge \mu'$,
and therefore the number of vertices in the first $j$~columns
but \itc{not} in the first $i$~rows
of the Young diagram of~$\lambda$ is at most
the total number of vertices in
rows $i+1$ through $i+j$ of the Young diagram of~$\lambda$.
Therefore $|C|$ is at most the sum of the first $i+j = k$ parts
of~$\lambda$.
\end{proof}

\begin{theorem}
\label{th:deltaconj}
If $\lambda$ is wide then the set of rows of the Young diagram of~$\lambda$
is a completely saturated clique cover of~$L(G_\lambda)$.
Moreover, $L(G_\lambda)$ satisfies $\Delta$-conjugacy
and the t-phenomenon.
\end{theorem}

\begin{proof}
By Propositions \ref{pr:uniform} and~\ref{pr:partition},
any uniform clique cover of the line graph of a bipartite graph
is completely saturated.
So in the case at hand,
Lemma~\ref{le:uniform} implies that the set of rows is completely
saturated.
The existence of a uniform clique cover implies that
$\Delta\omega$ is a partition, so the remaining claims
follow from Propositions \ref{pr:Saks} and~\ref{pr:partition}.
\end{proof}

The obvious remaining question is whether there exists a
uniform (or equivalently,
by Proposition~\ref{pr:uniform} and Theorem~\ref{th:deltaconj},
a completely saturated) stable set cover
of $L(G_\lambda)$ if $\lambda$ is wide. 
It is easy to see that the existence of such a cover is
equivalent to the WPC for free matroids.

\section{Relationship with network flows and with
edge colorings of bipartite graphs}

In the introduction we mentioned the existence of
a relationship between the WPC and
integer multicommodity flows (a.k.a.\ ``integer multiflows'').
To see this,
direct the edges of~$G_\lambda$ so that rows point to columns,
and give each edge a capacity of~one.
Enlarge $G_\lambda$ to a directed graph~$H_\lambda$
by adjoining $\lambda_1$ \itc{source vertices} $s_1, \ldots, s_{\lambda_1}$
and $\lambda_1$ \itc{destination vertices} $d_1, \ldots, d_{\lambda_1}$,
and adding a directed edge of capacity one from
each $s_i$ to each row of~$\lambda$
and from each column of~$\lambda$ to each~$d_i$.
What we seek is a simultaneous routing of $\lambda_1$ commodities
on~$H_\lambda$; specifically, we want to send
$\lambda'_i$ units of commodity~$i$ from $s_i$ to~$d_i$,
where the amount of every commodity on every link is
required to be an integer.

In this language, the WPC for free matroids essentially states that
the multiflow problems coming from wide partitions enjoy
a max-flow min-cut
property, and have integral optimal solutions.
Multiflow problems in general do not satisfy max-flow min-cut;
this is another way of seeing why the WPC for free matroids
cannot be proved purely by ``general nonsense,''
and that something is special about wide partitions
(if the conjecture is true).

The game of finding technical conditions to ensure max-flow min-cut
has been played before in the literature.
Unfortunately, we have been unable to find anything
that applies directly to our situation;
the graph $H_\lambda$ does not satisfy
any kind of Eulerian condition or topological condition
that is known to be helpful.
It is also readily seen that
the coefficient matrix of the linear programming relaxation
of this multiflow viewpoint is not totally unimodular.

Nevertheless, we are able to obtain some partial results,
which we present now.

\begin{lemma}
\label{Wide}
A partition $\lambda$ is wide if and only if for $L(G_{\lambda})$,
$$ \Delta \alpha = \lambda'.$$
\end{lemma}

\begin{proof}
The ``only if'' part follows from the results of the previous section,
but we ignore this and give a self-contained proof.
We show that being wide is equivalent to the condition
$$ \forall k: \ \ \alpha_k = \sum_{j=1}^{k}{\lambda_j'}. $$

We construct a directed network by taking $G_\lambda$ with edges directed from
the row vertices to the column vertices and with capacity $1$,
adding a source $s$ connected to each row vertex by an edge of capacity
$k$ and a target $t$ connected from each column vertex by an edge of capacity
$k$. The maximum flow in this network has value exactly $\alpha_k$, because
$k$-stable sets in the line graph correspond to edge subsets of
$G_\lambda$ of maximum
degree $k$ (since line graphs of bipartite graphs are perfect).

Consider a cut $C=(S,S')$ in this network ($s \in S, t \in S'$).
First choose $R$, the row vertices in $S$. The optimal way to add column
vertices to $S$ is to include $y \in S$ if it has at least $k$ neighbors
in $R$ (because then it is cheaper to have the edge $(y,t)$ in the cut rather
than the edges from $y$'s neighbors in $R$ to it). Thus the weight of the
minimum cut $C_R$ for a given $R$ is
$$ w(C_R) = k (n - |R|) + \sum_{j}{\min\{k, |N(j) \cap R|\}} $$
where $n$ is the number of rows and $N(j)$ is the set of neighbors
of column vertex $j$. $|N(j) \cap R|$ is the size of the $j$-th column
of the subpartition defined by $R$.

If the partition is wide, we have $\sum_{j}{\min\{k, |N(j) \cap R|\}} \geq
\sum_{j=1}^{k}{|N(j) \cap R|}$ and thus
$$ w(C_R) \geq k (n - |R|) + \sum_{j=1}^{k}{|N(j) \cap R|} \geq
\sum_{j=1}^{k}{\lambda_j'} $$
which means that the minimum cut is at least $\sum_{j=1}^{k}{\lambda_j'}$.
On the other hand, this value is achieved by setting $S$ to contain the
vertices corresponding to rows of length at most $k$. By the
max-flow min-cut theorem, the maximum flow is equal to
$\sum_{j=1}^{k}{\lambda_j'}$.

Conversely, if $\alpha_k = \sum_{j=1}^{k}{\lambda_j'}$, consider
a $k$-stable set $F_k$ of size $\alpha_k$. Since any $k$-stable set
has at most $\min\{k, \lambda_i\}$ squares in each row $i$ and
$\alpha_k=\sum_i \min\{k, \lambda_i\}=  \sum_{j=1}^{k}{\lambda_j'}$,
we have that $F_k$ has exactly $\min\{k, \lambda_i\}$ squares in each
row $i$. Consider now any subset of rows $R$, and
let $G_k$ be the restriction of $F_k$ to rows $R$. Thus the size
of $G_k$ is the size of the first $k$ columns intersected with $R$.
On the other hand, $G_k$ has at most $k$ squares in each column,
therefore its size is at most that of the first $k$ rows of $R$.
Over all $k$ and all subsets $R$, this implies that $\lambda$ is wide.
\end{proof}

\begin{lemma}
\label{Chains}
Let $G$ be the line graph of a bipartite graph,
and let $b$ be the number of distinct part sizes of $\Delta \alpha(G)$.
Let $a_1 > a_2 > \cdots > a_b$
be these part sizes and $k_i$ the number of parts of size $\geq a_i$.
Then a uniform stable set cover exists if and only if there exists
a chain
$$ F_1 \subset F_2 \subset \cdots \subset F_b $$
where $F_i$ is a $k_i$-stable set of size $\alpha_{k_i}$.
\end{lemma}

\begin{proof}
It is easy to see that if $(A_1, A_2, \ldots, A_{k_b})$ is a uniform
stable set cover, then
$$ F_i = \bigcup_{k=1}^{k_i}{A_k} $$
is a $k_i$-stable set of size $\alpha_{k_i}$ and these sets form a chain.

Conversely, suppose that we have such a chain $\emptyset = F_0 \subset
F_1 \subset \cdots \subset F_b$. Now consider each $F_i$ as a set of
edges in the underlying bipartite graph. Define $g_i$ to be the maximum
degree in $G_i = F_i \setminus F_{i-1}$. We would like to have
$g_i \leq k_i - k_{i-1}$ for each $i$. Therefore, take a chain where
the vector $(g_1, g_2, \ldots, g_b)$ is lexicographically minimal
and assume that $j$ is the first index where $g_j > k_j - k_{j-1}$.
Note that $\forall i < j: g_i = k_i - k_{i-1}$, otherwise $F_{j-1}$
would have degrees strictly smaller than $k_{j-1}$. Then it could be
extended to a larger $k_{j-1}$-stable set in the line graph. But $F_{j-1}$
is by assumption the maximum $k_{j-1}$-stable set. Also, $G_1 = F_1$
has degrees at most $k_1$, therefore $g_1 = k_1$ and $j > 1$.

Let $x$ be a vertex with degree $g_j$ in $G_j$. Since $g_j > k_j - k_{j-1}$
and $F_j$ has degrees at most $k_j$, $x$ has degree strictly smaller than
$k_{j-1}$ in $F_{j-1}$. Assume $x$ is on the ``left-hand side".
Consider all paths from $x$, using edges from $G_j$ and $G_{j-1}$
alternately. Let $H$ denote the union of all these paths. We claim
that for any vertex $y$ on the right-hand side, reachable from $x$ in $H$,
\begin{itemize}
\item $y$ has degree $\geq k_{j-1} - k_{j-2}$ in $G_{j-1}$.
\item $y$ has degree $\leq k_j - k_{j-1}$ in $G_j$.
\end{itemize}
By contradiction, if either of these conditions were violated, $y$ would
have degree strictly smaller than $k_{j-1}$ in $F_{j-1}$. (This follows
from the assumptions on $F_{j-2}$ and $F_j$.) Then we could switch the
edges on the (odd length) $x-y$ path between $G_{j-1}$ and $G_j$, thereby
increasing the size of $F_{j-1}$, while it would remain a $k_{j-1}$-stable
set in the line graph. However, $F_{j-1}$ had size $\alpha_{k_{j-1}}$
which was maximum.

This implies that we can estimate the number of edges in $G_j \cap H$
and $G_{j-1} \cap H$. The degrees in $G_{j-1}$ on the right are actually
equal to $k_{j-1} - k_{j-2}$, because $j$ is the first index where a
higher degree exists. Thus if there are $r$ vertices on the right-hand
side, reachable in $H$, we have
$$ |G_{j-1} \cap H| = r (k_{j-1} - k_{j-2}), $$
$$ |G_j \cap H| \leq r(k_j - k_{j-1}).$$

However, there is a vertex on the left-hand side ($x$) which has degree
strictly greater than $k_j - k_{j-1}$ in $G_j \cap H$. By assumption,
every vertex on the left has degree at most $k_{j-1} - k_{j-2}$
in $G_{j-1}$, so there must be a vertex $z$ on the left, reachable in $H$,
which has degree strictly smaller than $k_j - k_{j-1}$ in $G_j \cap H$.
By switching the edges between $G_j$ and $G_{j-1}$ on the path from $x$
to $z$, we maintain all the properties of $F_{j-1}$ and $F_j$; however,
we have decreased the number of vertices of degree $g_j$ in $G_j$. If
there are still vertices of degree $g_j$ in $G_j$, we repeat this procedure
until we decrease the maximum degree to $g_j - 1$. For each $i < j$,
we have maintained $g_i = k_i - k_{i-1}$. This contradicts the assumption
that the vector $(g_1, g_2, \ldots, g_b)$ is lexicographically minimal.

Now we have a chain $\emptyset = F_0 \subset F_1 \subset \cdots
\subset F_b$ where the degrees in $G_i = F_i \setminus F_{i-1}$
are at most $k_i - k_{i-1}$.
By Birkhoff-von Neumann, we can decompose each $G_i$ into
$k_i - k_{i-1}$ matchings $A_i^{(1)}, A_i^{(2)}, \ldots,
A_i^{(k_i-k_{i-1})}$. Each of these matchings must have size $a_i$
otherwise the largest one together with $F_{i-1}$ would form
a $(k_{i-1}+1)$-stable set larger than $\alpha_{k_{i-1}} + a_i =
\alpha_{k_{i-1}+1}$.

We have constructed a stable set cover
$$ A_1^{(1)}, A_1^{(2)}, \ldots, A_2^{(1)}, A_2^{(2)}, \ldots,
A_b^{(1)}, \ldots, A_b^{(k_b-k_{b-1})} $$
where the sizes of the stable sets are exactly the parts of $\Delta \alpha$.
\end{proof}

To see the power of the above lemmas,
first note that Proposition~\ref{pr:weak} follows easily.

\begin{proof}[Alternative proof of Proposition~\ref{pr:weak}]
Consider a weakly Latin tableau. Define $F_k$ to be the set
of all squares containing numbers up to $k$. Now consider $F_k$ as
a set of edges in the bipartite graph. Since the degrees in $F_k$ are
at most $k$, it can be decomposed into $k$ matchings and therefore
$F_k$ is a $k$-stable set in the line graph. The size of $F_k$ is
$\sum_{j=1}^{k}{\lambda_j'}$ which is the maximum possible size
of a $k$-stable set. By Lemma~\ref{Chains}, there exists a uniform stable
set cover, which corresponds to a Latin tableau.
\end{proof}

We can also easily deduce the following result.

\begin{theorem}
\label{th:twoparts}
If $\lambda$ is a wide partition with at most two distinct part sizes,
then $\lambda$ is Latin.
\end{theorem}

\begin{proof}
Let a partition $\lambda$ have parts of two different sizes $k_1 < k_2$.
By Lemma~\ref{Wide}, $\Delta \alpha = \lambda'$ which has $k_1$ parts of
one size and $k_2 - k_1$ parts of another (smaller) size.
There is a $k_1$-stable set of size $\alpha_{k_1}$ and
a $k_2$-stable set of size $\alpha_{k_2}$. The latter is the complete
set of vertices, so they form a chain trivially. By Lemma~\ref{Chains},
there exists a uniform stable set cover, which corresponds to a Latin
tableau.
\end{proof}

It is worth mentioning that Theorem~\ref{th:twoparts}
also follows from known results on edge colorings of bipartite graphs,
in particular from the following result of Folkman and Fulkerson.

\begin{definition}
Let $A$ be an $m\times n$ 0-1 matrix with a total of $N$~1's.
Let $\mu$ be a partition of~$N$.
We say that $A$ is \itc{$\mu$-decomposable} if
$A$ can be written as a sum
$$A = P_1 + P_2 + \cdots + P_{\ell(\mu)}$$
of 0-1 matrices~$P_i$ such that for all~$i$,
$P_i$ has a total of exactly $\mu_i$ 1's
and has at most one~1 in each row and column.
\end{definition}

\begin{proposition}[Folkman and Fulkerson]
\label{pr:FF}
Let $A$ be an $m\times n$ 0-1 matrix with a total of $N$~1's.
Let $\mu$ be a partition of~$N$
with at most two distinct part sizes.
Then $A$ is $\mu$-decomposable
if and only if
every $e\times f$ submatrix $B$ of~$A$ has at least
the following number of 1's:
$$\sum_{i\ge (m-e) + (n-f) + 1} \mu'_i.$$
\end{proposition}

\begin{proof}
Theorem 3.1 of \cite{FF}.
\end{proof}

\begin{proof}[Alternative proof of Theorem~\ref{th:twoparts}]
Let $m=\ell(\lambda)$ and let $n=\lambda_1$.
Let $A$ be the $m\times n$ matrix whose $(i,j)$ entry is~1
if $(i,j)$ is a cell of~$\lambda$
(i.e., if $j\le\lambda_i$)
and whose $(i,j)$ entry is~0 otherwise.
Let $\mu=\lambda'$.
Then $\mu$ also has at most two distinct part sizes.
Chasing definitions, we see that $A$ is $\mu$-decomposable
if and only if $\lambda$ is Latin.
We therefore need only check that the wideness of~$\lambda$
implies that the condition on submatrices of~$A$ in Proposition~\ref{pr:FF}
is satisfied. 
This is straightforward and we leave the details to the reader.
\end{proof}

It is tempting to wonder how far
Theorem~\ref{th:twoparts} may be generalized.
Perhaps the WPC for free matroids can be generalized to
arbitrary bipartite graphs?
Unfortunately, the answer is no;
if the condition on the number of distinct part sizes of~$\mu$
in Proposition~\ref{pr:FF} is dropped,
then it no longer remains true, and
a counterexample may be found in~\cite{FF}.
However, it is possible that as far as edge colorings are concerned,
it is being a partition 
that is the crucial property (rather than being wide).
More precisely, the following question remains open.

{\bf Latin Tableau Question, alternative form.}
Does Proposition~\ref{pr:FF} remain true if
the condition on the number of distinct part sizes of~$\mu$ is dropped
but $A$ is required to arise from a Young diagram (i.e., $A$ must
satisfy the condition that whenever $A_{ij} = 1$ then $A_{rs} = 1$
for all $r\le i$ and $s\le j$)?

It is not hard to show that this question is indeed equivalent to
the Latin Tableau Question as previously formulated.
Surprisingly, in spite of the sizable literature on
edge colorings of bipartite graphs,
the condition that $A$ arise from a Young diagram
does not seem to have been directly addressed before.

The set of all color-feasible partitions
(i.e., partitions~$\mu$ for which there exists an edge coloring
in which color~$i$ is used exactly $\mu_i$ times)
for a given bipartite graph
does not in general have a unique maximal element in dominance order.
But as we mentioned in the introduction,
the WPC for free matroids is equivalent to the claim
that for~$G_\lambda$ (with $\lambda$ wide),
there \itc{is} a unique maximal element.
Now, a necessary and sufficient condition
for the existence of a unique maximal element is given in~\cite{deWerra}.
Unfortunately, this necessary and sufficient condition
does not seem easy to verify for wide partitions.
However, the main theorem of~\cite{deWerra} does imply the following.

\begin{theorem}
If $\lambda$ is a wide partition with three distinct part sizes
and either the second or third part size occurs with multiplicity one, or
if $\lambda$ is a wide partition with four distinct part sizes
and the second and fourth part sizes both occur with multiplicity one,
then $\lambda$ is Latin.
\end{theorem}

\begin{proof}
This may be deduced from Corollary~3.3 of~\cite{deWerra}
in the same manner that
we deduced Theorem~\ref{th:twoparts} from Proposition~\ref{pr:FF}.
\end{proof}

We have one final result along the same lines.

\begin{theorem}
\label{th:threeparts}
If $\lambda$ is a self-conjugate wide partition
with at most three distinct part sizes,
then $\lambda$ is Latin.
\end{theorem}

\begin{proof}
Let $\lambda$ be a self-conjugate wide partition
with exactly three distinct part sizes.
(The case of one part size is trivial and
the case of two part sizes is covered by Theorem~\ref{th:twoparts}.)
Let $m_1$ be the multiplicity of the largest part size,
let $m_2$ be the multiplicity of the next largest part size, and
let $m_3$ be the multiplicity of the smallest part size.
Call the integers from $1$ to~$m_1$ the \itc{low range,}
call the integers from $m_1+1$ to $m_1+m_2$ the \itc{mid range,} and
call the integers from $m_1+m_2+1$ to $m_1+m_2+m_3$ the \itc{high range}.

The Young diagram of~$\lambda$ subdivides naturally into
six rectangular subregions,
which we give names as shown in the picture below.

\def\pstrut{$\phantom{\biggl|}$}
$$\vbox{\tabskip=0pt \offinterlineskip
 \halign to 15em{\tabskip=0pt plus1em#&#&\hfil#\hfil
   &#&\hfil#\hfil&#&\hfil#\hfil&#\tabskip=0pt\cr
 &\multispan7\hrulefill\cr
 &\vrule&\pstrut $A$&\vrule&$B$&\vrule&$D$&\vrule\cr
 &\multispan7\hrulefill\cr
 &\vrule&\pstrut $B'$&\vrule&$C$&\vrule\cr
 &\multispan5\hrulefill\cr
 &\vrule&\pstrut $D'$&\vrule\cr
 &\multispan3\hrulefill\cr
}
}$$
In addition, we define $E$ to be the square region $A \cup B\cup B'\cup C$.

In view of Lemma~\ref{Chains},
it suffices to construct
a subset $\alpha\subseteq A$
containing exactly $m_3$ cells from each row and each column of~$A$,
and a subset $\beta \subseteq E$,
disjoint from~$\alpha$, containing
exactly $m_2$ cells from each row and each column of~$E$.
We split into two cases.

{\em Case 1:} $m_1 \ge m_2+m_3$.
Temporarily place any $m_1\times m_1$ Latin square~$L$ into region~$A$.
(The only purpose of~$L$ is to help describe $\alpha$ and~$\beta$.)
Let $\alpha$ be the set of cells of~$L$ with an entry between
1 and~$m_3$ inclusive.
Let $b$ be the set of cells of~$L$ with an entry between
$m_3+1$ and~$m_3+m_2$ inclusive,
and let $\beta = b\cup C$.
It is easily checked that $\alpha$ and~$\beta$ have the desired properties.

{\em Case 2:} $m_1 < m_2+m_3$.
The set~$\alpha$ may be constructed exactly as in Case~1,
but the construction of~$\beta$ requires several steps.

Let $b$ be a subset of~$B$ with the following two properties:
$(1)$~each row of $b$ contains $m_2+m_3-m_1$ cells, and
$(2)$~the number of cells in
any two columns of~$b$ differ by at most one.
It easy to see that the Gale-Ryser theorem implies that
such a subset~$b$ exists.

Let $c_i$ be the number of cells in the $i$th column of~$b$.
We claim that $c_i\le m_2$ for all~$i$.
To see this, note that $\sum_i c_i = m_1(m_2+m_3-m_1)$.
Since any two $c_i$ differ by at most one,
it follows that if $c_i>m_2$ for some~$i$ then $c_j\ge m_2$ for all~$j$.
Since $B$ has $m_2$~columns,
it follows that $\sum_i c_i > m_2^2$.
Therefore,
$m_2^2 < m_1(m_2+m_3-m_1)$.
However, we claim that the wideness of $\lambda$ implies that
\begin{eqnarray}
\label{eq:star}
m_1^2 + m_2^2 \ge m_1(m_2+m_3),
\end{eqnarray}
yielding the desired contradiction.
To see why the inequality~(\ref{eq:star}) is true,
suppose first that $m_2 < m_1$.
The lower subpartition $B' \cup C\cup D'$ of~$\lambda$ is wide,
so in particular the sum of its first $m_1$ rows
is at least the sum of its first $m_1$ columns.
Then inequality~(\ref{eq:star}) follows immediately.
On the other hand, suppose $m_2\ge m_1$.
The rectangle~$D'$ is wide, so $m_1\ge m_3$.  Therefore,
$$m_1^2 + m_2^2 \ge m_1^2 + m_1m_2 = m_1(m_1+m_2) \ge m_1(m_2+m_3),$$
yielding inequality~(\ref{eq:star}) again.

Since $c_i \le m_2$, the quantity $m_2 - c_i$ is
a nonnegative integer for all~$i$.
Since any two $c_i$ differ by at most one,
another easy application of Gale-Ryser implies that
there exists a subset $c \subseteq C$ whose
$i$th row contains exactly $m_2 - c_i$ cells and whose
$i$th column also contains exactly $m_2 - c_i$ cells.

Finally, we set
$$\beta = (A\backslash \alpha) \cup b \cup b' \cup c,$$
where $b'$ is the subset of~$B'$ that is the transpose of~$b$.
Again one easily checks that $\alpha$ and~$\beta$ have the required
properties.
\end{proof}

\section{Reduction to self-conjugate partitions}

\begin{theorem}
\label{SelfConj}
Let $\lambda = (\lambda_1, \ldots, \lambda_n)$ be a wide partition,
and let $m = \lambda_1$.
Let $\mu$ be the following partition with $2m + n$ parts:
$$ \mu = (2m + \lambda_1', \ldots, 2m + \lambda_m', m, \ldots, m,
\lambda_1, \ldots \lambda_n).$$
(In other words, $\mu$ is a $2m \times 2m$ square with $\lambda$ added
on the bottom and $\lambda'$ added on the right.)
Then $\mu$ is a wide partition.
\end{theorem}

\begin{proof}
We use Lemma~\ref{Wide} and prove that for any $k$, there is a $k$-stable
set in $L(G_{\mu})$ of size $\sum_{j=1}^{k}{\mu_j'}$. We distinguish
three cases:

{\em Case 1: $k \leq m$.}
We know $L(G_{\lambda})$ has a $k$-stable set of size $\sum_{j=1}^{k}
{\lambda_j'}$. Denote this set by $F$. We define a $k$-stable set $F'$
in $L(G_\mu)$:
First, include $(2m+i,j) \in F'$ and $(j,2m+i) \in F'$ for each $(i,j)
\in F$. To define the remaining part of $F'$ (in the $2m \times 2m$
square), we need to find a bipartite graph on $2m + 2m$ vertices with
a given sequence of degrees on both sides: $m$ degrees equal to $k$
and the remaining degrees smaller than $k$.

\begin{figure}[ht]
\label{Case1}
\caption{}
\begin{center}
\scalebox{0.4}{\includegraphics{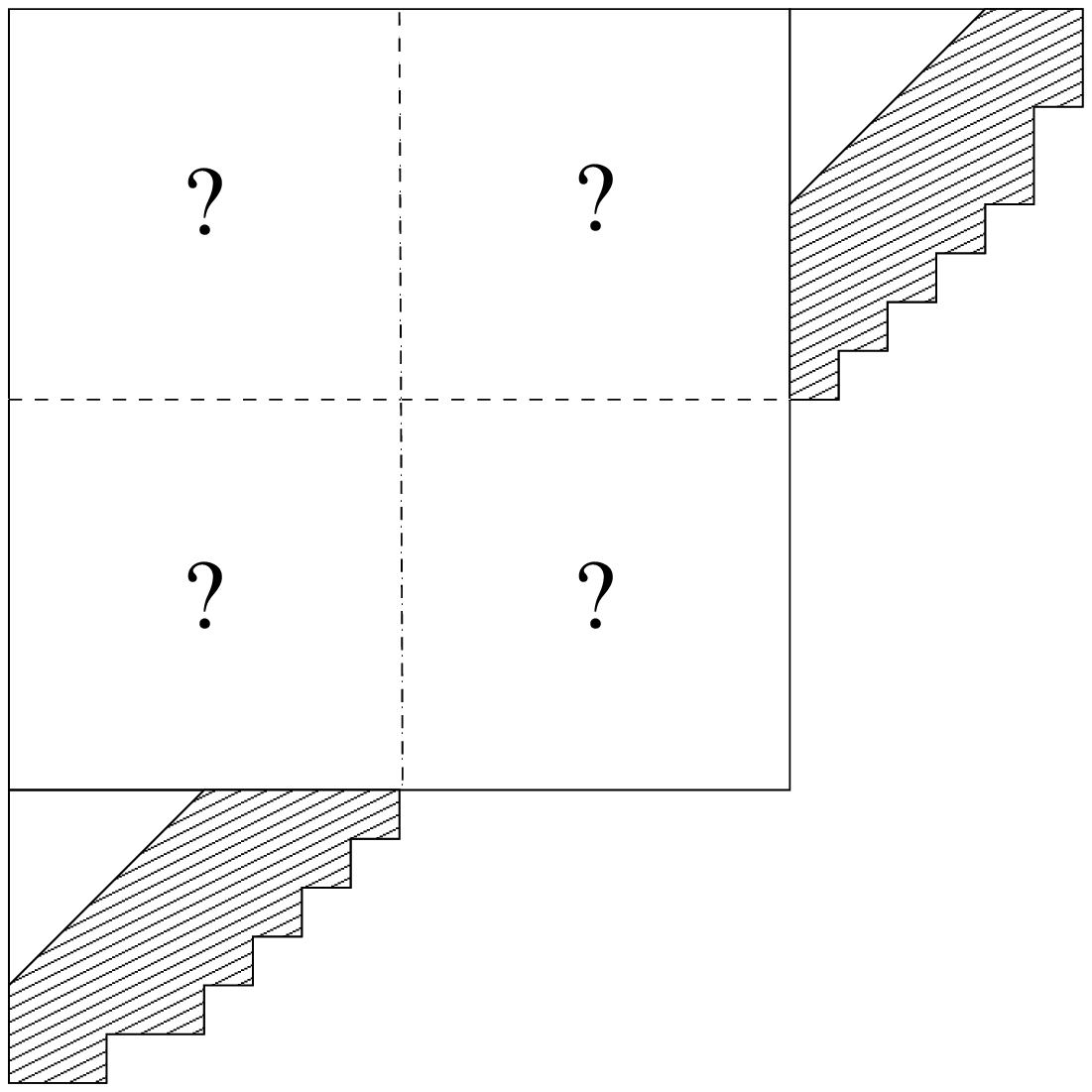}}
\end{center}
\end{figure}

We find the bipartite graph using the Gale-Ryser theorem
(Proposition~\ref{pr:GR}), which may be restated as follows.
There is a bipartite graph with degrees $\sigma_1 \geq \sigma_2 \geq
\cdots \geq \sigma_p$ on the left and $\rho_1 \geq \rho_2 \geq \cdots
\geq \rho_p$ on the right, if and only if $\sigma$ and $\rho$ as
partitions satisfy
$$ \sigma' \geq \rho.$$

In this case, we have $\sigma=\rho$ and $\sigma_1 = \cdots = \sigma_m
= k$, i.e.  $\forall i: \sigma_i' \geq m \geq k$. On the other hand,
$\forall i: \rho_i \leq k$ which implies that $\sigma' \geq \rho$.

{\em Case 2: $m < k \leq 2m$.}
In this case, we include in $F'$ all squares $(i,j)$ with either
$i > 2m$ or $j > 2m$. Also, we include the squares $(m+i,m+j)$ for
$1 \leq i, j \leq m$ and squares $(m+i,j)$ and $(j,m+i)$ satisfying
$(j-i) \bmod m \in\{0,1,\ldots,k-m-1\}$. 
To complete $F'$, we must find a bipartite graph
on $m + m$ vertices (the top-left $m \times m$ square) with degrees
on both sides equal to $d_i = m - \lambda_i'$.

\begin{figure}[ht]
\label{Case2}
\caption{}
\begin{center}
\scalebox{0.4}{\includegraphics{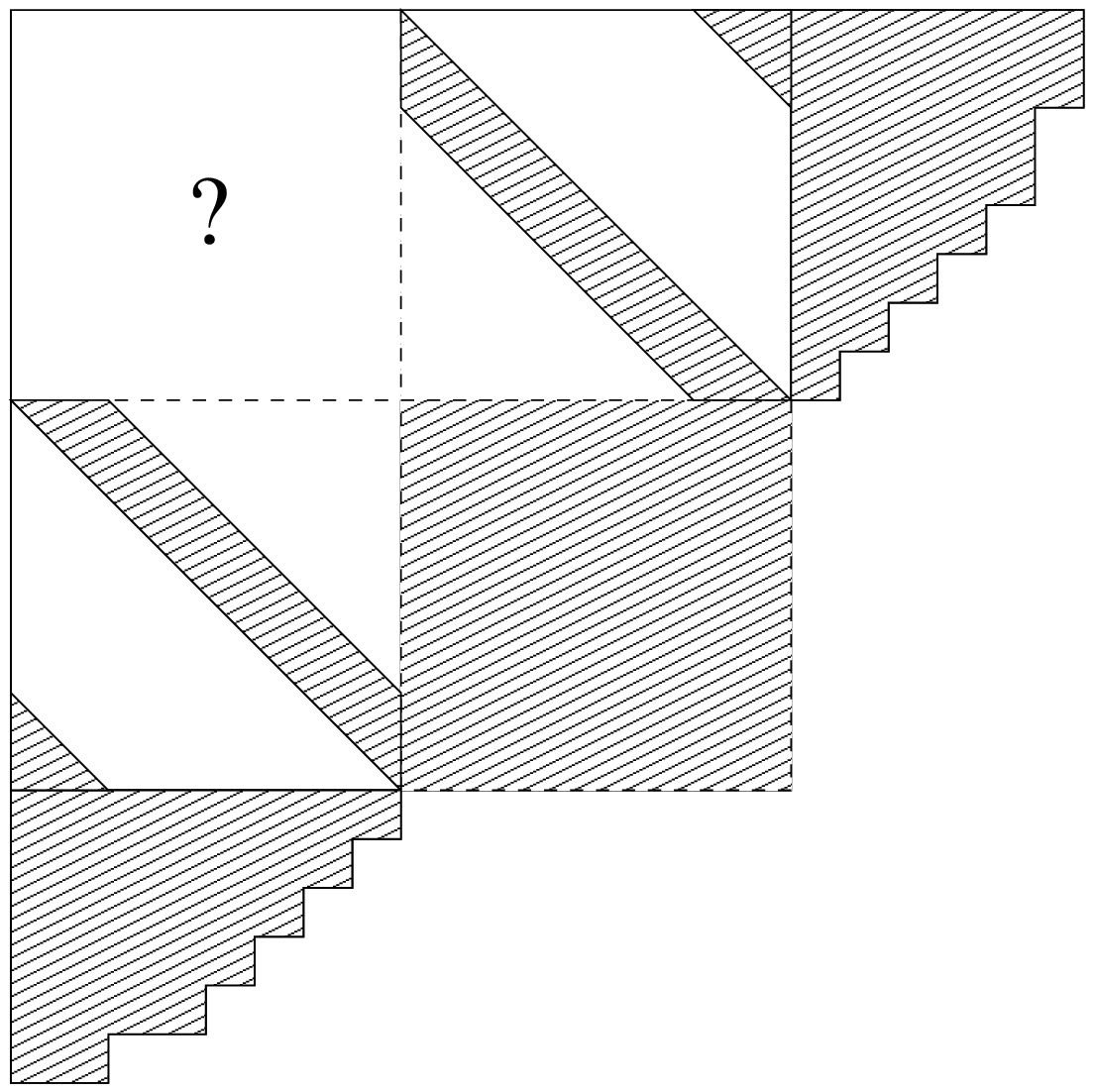}}
\end{center}
\end{figure}

Again, we apply the Gale-Ryser theorem. We find the complement of the
required bipartite graph, which should have degrees $m - d_i = \lambda_i'$
on both sides. Here $\sigma = \rho = \lambda'$ and $\lambda \geq \lambda'$
because $\lambda$ is a wide partition.

{\em Case 3: $k > 2m$.}
Here, we include all squares $(i,j)$ with $i>m$ or $j>m$.
To complete $F'$, we must find a bipartite graph on $m + m$ vertices
with degrees on both sides equal to $d_i = \min\{m, k - m - \lambda_i'\}$.

\begin{figure}[ht]
\label{Case3}
\caption{}
\begin{center}
\scalebox{0.4}{\includegraphics{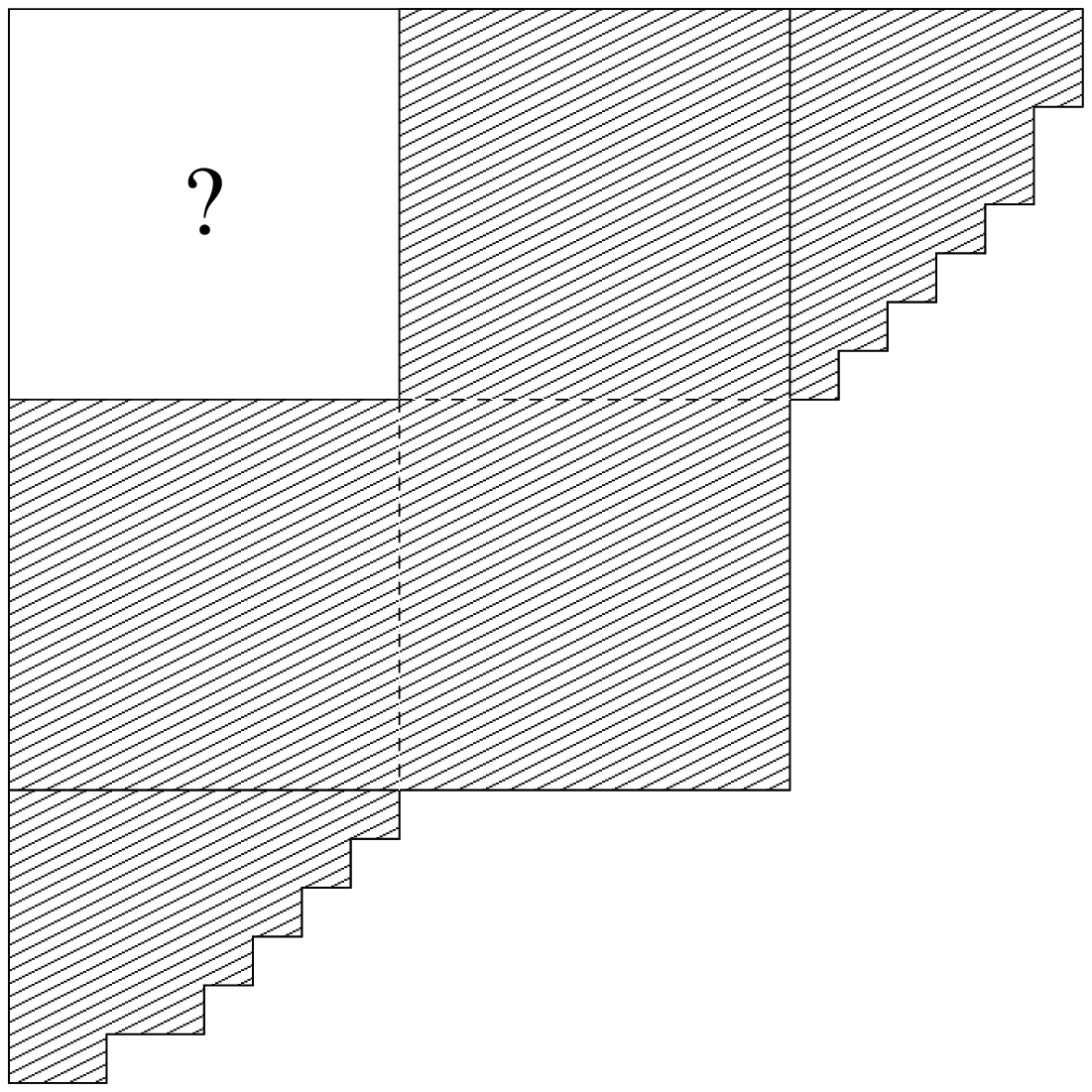}}
\end{center}
\end{figure}

Similarly to Case 2, we find the complement of the bipartite graph which
should have degrees $m - d_i = \max\{\lambda_i' - (k-2m), 0\}$ on both
sides. Here $\sigma' = \rho'$ is equal to $\lambda$ without the first
$k-2m$ rows. Since $\lambda$ is wide, again $\sigma' \geq \rho$.
\end{proof}

\begin{corollary}
If the Wide Partition Conjecture holds for self-conjugate wide partitions,
then it is true for all wide partitions.
\end{corollary}

\begin{proof}
Let $M$ be a matroid, $\lambda$ a wide partition and $I_i$ an independent
set given for each row. We define a self-conjugate wide partition $\mu$
containing $\lambda$ as above. We assign the same set $I_i$ to each row
of $\lambda$. We assign arbitrary independent sets to the remaining rows.
(If necessary, we extend the matroid to a sufficiently large $M'$ such
that $A$ is independent in $M'$ iff $A \cap M$ is independent in $M$.)

Assume that the Wide Partition Conjecture holds for self-conjugate
partitions. Then there exists a permutation of $I_i$ in each row
so that the set in each column is independent. Obviously, the assignment
restricted to $\lambda$ satisfies the same property.
\end{proof}

\section{Counterexamples}

One might hope that even for wide partitions with more than two part
sizes, one could build the desired chain of $k$-stable sets greedily,
either from the top or from the bottom. However, this is impossible,
since some maximum $k_i$-stable sets cannot be extended to any maximum
$k_{i+1}$-stable set and some maximum $k_i$-stable sets do not contain
any maximum $k_{i-1}$-stable set.

\begin{figure}[ht]
\caption{}
\label{fig:maxflow}
\begin{center}
\scalebox{0.5}{\includegraphics{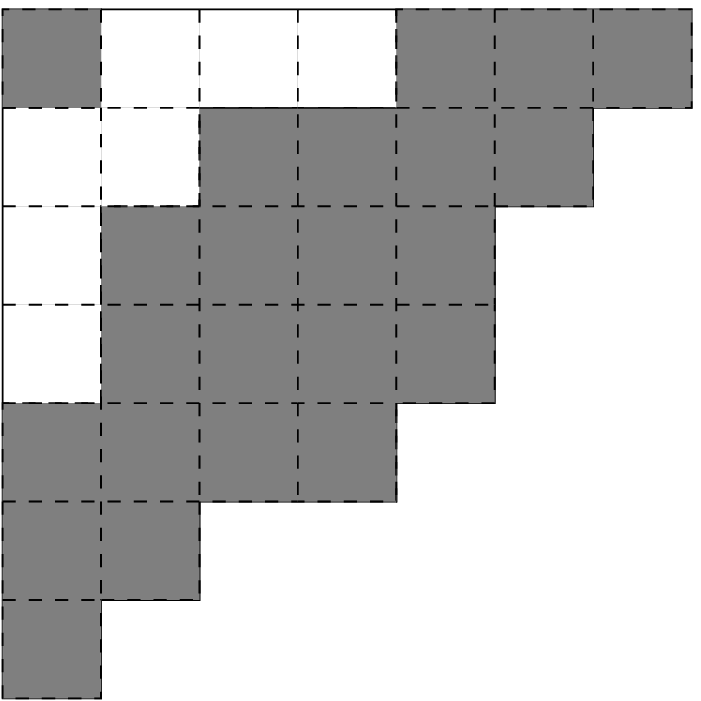}}
\end{center}
\end{figure}

Figure~\ref{fig:maxflow} shows
a maximum 4-stable set that is not extendible
to any maximum 5-stable~set.

\begin{figure}[ht]
\caption{}
\label{fig:minflow}
\begin{center}
\scalebox{0.5}{\includegraphics{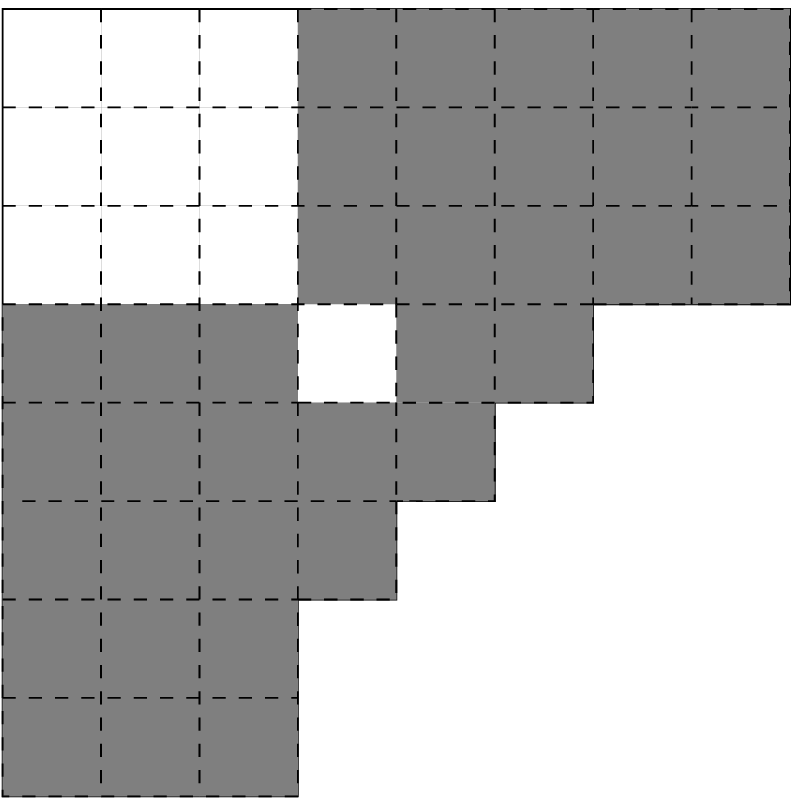}}
\end{center}
\end{figure}

Figure~\ref{fig:minflow} shows a maximum 5-stable set that contains
no maximum 4-stable set.

As we mentioned before,
uniform stable set covers do not always exist for line graphs of bipartite
graphs. Even for graphs of some ``skew shapes" (differences
of two partitions), there may be no chain of $k$-stable
sets along the lines of Lemma~\ref{Chains}.

\begin{figure}[ht]
\caption{}
\label{fig:twostable}
\begin{center}
\scalebox{0.5}{\includegraphics{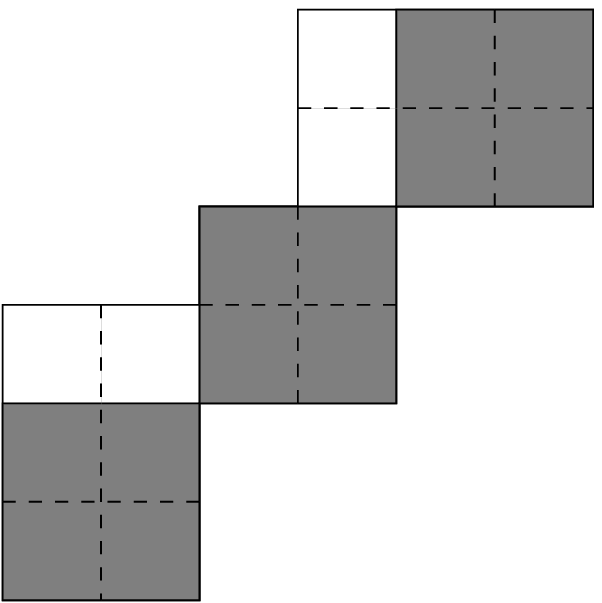}}
\end{center}
\end{figure}

\begin{figure}[ht]
\caption{}
\label{fig:threestable}
\begin{center}
\scalebox{0.5}{\includegraphics{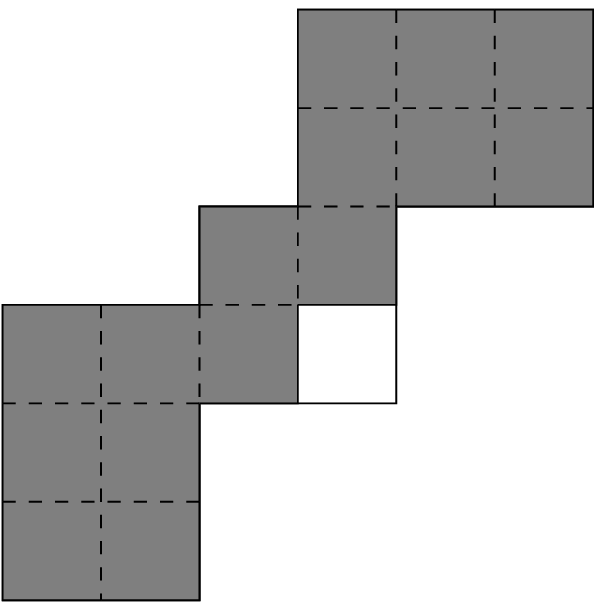}}
\end{center}
\end{figure}

For example,
the shaded area in Figure~\ref{fig:twostable}
is the unique maximum 2-stable set,
while the shaded area in Figure~\ref{fig:threestable}
is the unique maximum 3-stable set.
Thus there is no chain of maximum $k$-stable sets.

On a different note, it is tempting to try to prove
the WPC for free matroids
by explicitly filling in the Young diagram of~$\lambda$ one
row at a time or even one entry at a time.
Some such approach may indeed work,
but we have tried several such constructions without success.
For example, Sandy Kutin (personal communication)
has suggested filling in the rows one at a time
starting from the bottom, and whenever there is a choice, choosing the
lexicographically largest possibility.
This method fails for $\lambda=(6,6,6,5,2,2)$, as seen below.

\def\pstrut{$\phantom{\bigl|}$}
$$\vbox{\tabskip=0pt \offinterlineskip
 \halign to 10em{\tabskip=0pt plus1em#&#&\hfil#\hfil
   &#&\hfil#\hfil &#&\hfil#\hfil
   &#&\hfil#\hfil &#&\hfil#\hfil
   &#&\hfil#\hfil &#&\hfil#\hfil
   &#\tabskip=0pt\cr
 &\multispan{13}\hrulefill\cr
 &\vrule&\pstrut?&\vrule&?&\vrule&?&\vrule&?&\vrule&?&\vrule&?&\vrule& &\cr
 &\multispan{13}\hrulefill\cr
 &\vrule&\pstrut$4$&\vrule&$6$&\vrule&$5$&\vrule&$1$&\vrule&$3$&\vrule
    &$2$&\vrule\cr
 &\multispan{13}\hrulefill\cr
 &\vrule&\pstrut$6$&\vrule&$5$&\vrule&$4$&\vrule&$3$&\vrule&$2$&\vrule
    &$1$&\vrule\cr
 &\multispan{13}\hrulefill\cr
 &\vrule&\pstrut$5$&\vrule&$4$&\vrule&$3$&\vrule&$2$&\vrule&$1$&\vrule\cr
 &\multispan{11}\hrulefill\cr
 &\vrule&\pstrut$1$&\vrule&$2$&\vrule\cr
 &\multispan5\hrulefill\cr
 &\vrule&\pstrut$2$&\vrule&$1$&\vrule\cr
 &\multispan5\hrulefill\cr
}
}$$

\section{Acknowledgments}

The first author thanks Victor Reiner for
suggesting that
the WPC, network flows, and the Greene-Kleitman theorem might be related,
and thanks Debra Waugh for useful suggestions.
Sergey Fomin, Alexander Postnikov, and Mihai Ciucu
provided good feedback on earlier presentations of this work.

\bibliographystyle{plain}
\bibliography{refs}
\end{document}